\documentclass{ws-procs9x6}

\begin{document}

\title{$G$-STRUCTURES DEFINED ON
PSEUDO-RIEMANNIAN MANIFOLDS}

\author{IGNACIO SANCHEZ-RODRIGUEZ}

\address{Department of Geometry
and Topology,
University of Granada,\\
Granada, E-18071, SPAIN\\
E-mail: ignacios@ugr.es\\
www.ugr.es/local/ignacios}

\begin{abstract}
Concepts and techniques from the theory
of $G$-structures of higher order are
applied to the study of certain structures
(volume forms, conformal structures,
linear connections and projective structures)
defined on a pseudo-Riemannian manifold.
Several relationships between the
structures involved
have been investigated. The operations
allowed on $G$-structures,
such as intersection, inclusion,
reduction, extension and prolongation,
were used for it.
\end{abstract}


\bodymatter

\section{Introduction and motivation}
A \emph{differential structure}
of a manifold $M$
is a $C^{\infty }$
maximal atlas and, indeed,
the charts of the atlas make up
the \emph{primordial} structure.
The idea of a
\emph{geometrical structure}
can be realized by
the concept of $G$-structure
when choosing the
allowable meaningful
classes of charts.

General relativity
is a physical theory,
which is heavily based on
differential geometry.
The fundamental mathematical tools
used by this theory to explain
and to handle gravity are the
geometrical structures.
The space-time is described
by a 4-dimensional manifold
with a \emph{Lorentzian metric field},
and the theory put the matter on space-time,
being mainly represented
by curves in the manifold
or by the overall
stress-energy tensor.

The law of inertia in the space-time
is translated into
a \emph{projective structure}
on the manifold,
which is provided by
the geodesics of the metric
in keeping with the
equivalence principle.
Furthermore, the space-time in
general relativity is
a dynamical entity because
the metric field is subject
to the Einstein field equations,
which almost equate Ricci
curvature with stress-energy of matter.

Other main structures are
the \emph{volume form}, that is used
to get action functionals
by integration over the manifold,
and the \emph{Lorentzian conformal structure},
that gives an account of
light speed invariance.
Different approaches to gravity
try to \emph{separate} the geometry
into independent compounds
to promote the understanding
about physical interpretation
of geometric variables.

The theory of $G$-structures
of higher order is possibly
the more natural framework
for studying the interrelations
involved among the relevant structures.
In a pseudo-Riemannian
manifold there are
defined unambiguously the
following structures:
\emph{volume form,
conformal structure,
pseudo-Riemannian metric,
symmetric linear connection
and projective structure}.
Volume, conformal and metric
structures are
$G$-structures
of first order,
but each of them lead to
a \emph{prolonged}
second order structure.
Symmetric linear connection
and projective structure
are \emph{inherently}
$G$-structures of second order.
We will try to clarify this unified description.

\section{The bundle of $r$-frames}
A differentiable manifold $M$
is a set of points
with the property that we can cover it
with the charts of a $C^{\infty }$
$n$-dimensional maximal atlas
$\mathcal{A}$.
The \emph{bundle of $r$-frames}
$\mathcal{F}^{r}M$ is
a quotient set
over $\mathcal{A}$.
Every class-point, an
$r$-frame,
collect the charts with equal
origin of coordinates
which produce identical $r$-th
order Taylor series expansion
of functions
(see Refs.~\refcite{Koba1,IgTes}).

An $r$-frame is an $r$-jet at $0$
of inverses of charts of $M$;
two charts are in
the same $r$-jet
if they have the same
partial derivatives up to
$r$-th order at the same
origin of coordinates.
Every $\mathcal{F}^{r}M$
is naturally equipped with
a \emph{principal bundle} structure
with respect to the group
$\mathrm{G}^{r}_{n}$
of $r$-jets
at 0 of diffeomorphisms
of $\mathbb{R}^{n}$,
$j^{r}_{0}\phi $,
with $\phi(0)=0$.

The group of the
\emph{bundle of 1-frames}
is $\mathrm{GL}(n,\mathbb{R})\cong
\mathrm{G}^{1}_{n}$.
Its natural representation
on $\mathbb{R}^{n}$
gives an \emph{associated bundle}
coinciding with the
tangent bundle $TM$.
In the end, we identify
$\mathcal{F}^{1}M$ with the
\emph{linear frame bundle} $LM$.
Other representations of
$\mathrm{G}^{1}_{n}$
on subspaces of the
tensorial algebra over
$\mathbb{R}^{n}$
give associated bundles
whose sections are
the well-known tensor fields.

The \emph{bundle of 2-frames}
$\mathcal{F}^{2}M$
is somehow more complicated.
Every 2-frame
is characterized by a
\emph{torsion-free
transversal $n$-subspace}
$H_{l}
\subset T_{l}
LM$.
It happens that
the chart's first
partial derivatives
fix $l\in LM$
and the second
partial derivatives
give the 'inclination'
of that $n$-subspace.
The group
$\mathrm{G}^{2}_{n}$
is isomorphic to
$\mathrm{G}^{1}_{n}
\rtimes \mathrm{S}^{2}_{n}$,
a semidirect product, with
{$\mathrm{S}^{2}_{n}$}
the additive group
of symmetric bilinear maps
of $\mathbb{R}^{n}\times
\mathbb{R}^{n}$
into $\mathbb{R}^{n}$;
the multiplication rule is
$(a,s)(b,t):=(ab,b^{-1}s(b,b)+t)$,
for $a,b\in \mathrm{G}^{1}_{n}$,
$s,t\in \mathrm{S}^{2}_{n}$.

Let $\mathfrak{g}$ denote
the Lie algebra of $G\subset
\mathrm{G}^{1}_{n}$,
then the named \emph{first prolongation
of} $\mathfrak{g}$ is
defined by
$\mathfrak{g}_{1}:=
\mathrm{S}^{2}_{n}\cap L(
\mathbb{R}^{n},\mathfrak{g})$.
We obtain that
$G\rtimes
\mathfrak{g}_{1}$
is a subgroup of
$\mathrm{G}^{1}_{n}
\rtimes \mathrm{S}^{2}_{n}\cong
\mathrm{G}^{2}_{n}$. This will be used
in Sec.~\ref{SeProl}.

\section{First order $G$-structures}
An \emph{$r$-th order
$G$-structure} of $M$ is
a \emph{reduction}
(see \cite[p. 53]{KoNo1})
of $\mathcal{F}^{r}M$ to a subgroup
$G\subset \mathrm{G}^{r}_{n}$.
First order $G$-structures are
just called $G$-structures.
Let us see some of them.

Let's define a
\emph{volume on $M$} as a
first order $G$-structure,
with $G= \mathrm{SL}^{\smash[t]{\pm}}_{n}:=\{a
\in \mathrm{G}^{1}_{n}\colon
|{\mathrm det}(a)|=1\}$.
If $M$ is orientable,
a volume on $M$ has two components:
two $\mathrm{SL}(n,
\mathbb{R})$-structures
for two equal,
except sign, \emph{volume $n$-forms}.
For a general $M$, a volume
corresponds to an
\emph{odd type $n$-form}.

From \emph{bundle
theory}, \cite{KoNo1}
$\mathrm{SL}^{\pm}_{n}$-structures
are sections of
$\mathcal{V}M$, the
\emph{associated bundle} to
$LM$ and the action of
$\mathrm{G}^{1}_{n}$ on
$\mathrm{G}^{1}_{n}/
\mathrm{SL}^{\pm}_{n}$,
and they correspond to
$\mathrm{G}^{1}_{n}$-\emph{equivariant
functions} $f$
of $LM$ to $\mathrm{G}^{1}_{n}/
\mathrm{SL}^{\pm}_{n}$.
The isomorphisms $\mathrm{G}^{1}_{n}/
\mathrm{SL}^{\pm}_{n}
\simeq \mathrm{H}_{n}:=
\{ k\, I_{n}\ :\ k>0\}
\simeq \mathbb{R}^{+}$
allow to write
$f\colon LM\rightarrow
\mathbb{R}^{+}$;
the equivariance
condition is
$f(la)=
|\det a|^{\frac{-1}{n}}f(l)$,
for $l\in LM$, $a\in
\mathrm{G}^{1}_{n}$.

\begin{theorem}\label{ThVol}
We have the bijections:
\[\text{Volumes on }M
\ \longleftrightarrow \
\operatorname{Sec}
\mathcal{V}M
\ \longleftrightarrow \
C^{\infty}_{\text{equi}}
(LM,\mathbb{R}^{+})\]
\end{theorem}
Analogous bijective diagram
can be obtained for every
reduction of a principal bundle.
The Lie algebra of
$\mathrm{SL}^{\pm}_{n}$ is
$\mathfrak{sl}(n,\mathbb{R})$
and its first prolongation
is $\mathfrak{sl}(n,\mathbb{R})_{1}
=\{s\in \mathrm{S}^{2}_{n}
\colon \sum _{k}s^{k}_{ik}=0\}$;
it's a Lie algebra of
infinite type.

We define a
\emph{pseudo-Riemannian metric}
as an
$\mathrm{O}_{q,n\textrm{-}q}$-structure,
with
$\mathrm{O}_{q,n\textrm{-}q}:=\left\{ a
\in \mathrm{G}^{1}_{n}\colon
a^{t}\eta a=\eta:=
\left( \begin{smallmatrix}
-I_{q}&0\\
0&I_{n\textrm{-}q}
\end{smallmatrix}\right) \right\} $.
As in Th.\ref{ThVol}, we
obtain bijections
between the metrics
and the sections of the
associated bundle
with typical fiber
$\mathrm{G}^{1}_{n}/
\mathrm{O}_{q,n\textrm{-}q}$,
and also with the
equivariant functions of $LM$ in
$\mathrm{G}^{1}_{n}/
\mathrm{O}_{q,n\textrm{-}q}$.
The first prolongation of
$\mathfrak{o}_{q,n\textrm{-}q}$
is $\mathfrak{o}_{q,n\textrm{-}q\ 1}=0$;
a consequence of this fact is
the uniqueness of the \emph{Levi-Civita
connection}.

A (pseudo-Riemannian)
\emph{conformal structure} is a
$\mathrm{CO}_{q,n\textrm{-}q}$-structure,
with $\mathrm{CO}_{q,n\textrm{-}q}
:=\mathrm{O}_{q,n\textrm{-}q}
\cdot \mathrm{H}_{n}$ (direct product);
this definition is equivalent
to consider a class of metrics
related by a positive factor,
and in the Lorentzian case, $q=1$,
a conformal structure
is characterized by
the \emph{field of null cones}.
The first prolongation
of $\mathfrak{co}_{q,n\textrm{-}q}$
is $\mathfrak{co}_{q,n\textrm{-}q\ 1}=
\{ s\in \mathrm{S}^{2}_{n}\colon
s^{i}_{jk}=\delta^{i}_{j}\mu _{k}
+\delta^{i}_{k}\mu _{j}
-\sum_{s}\eta^{is}\eta_{jk}\mu _{s},\
\mu =(\mu _{i})\in
\mathbb{R}^{n\,\ast}\} \simeq
\mathbb{R}^{n\,\ast}$.
The named second prolongation
$\mathfrak{co}_{q,n\textrm{-}q\ 2}$
is equal to 0 (i. e.,
$\mathfrak{co}_{q,n\textrm{-}q}$
is of finite type 2); this
deals with the existence and uniqueness of
the normal Cartan connection
but we do not deal with this here
(see \cite[\S\S VI.4.2, VII.3]{IgTes}).

Volumes on $M$
and conformal structures
are \emph{extensions}
(see \cite[p. 202]{GHV1})
of pseudo-Riemannian metrics
because of the inclusion
of $\mathrm{O}_{
q,n\textrm{-}q}$ in
$\mathrm{SL}^{
\pm}_{ n}$
and $\mathrm{CO}_{
q,n\textrm{-}q}$.
Reciprocally:

\begin{theorem}\label{ThMetric}
\ A pseudo-Riemannian metric
field on $M$
is  given by
a pseudo-Riemannian
conformal structure
and a volume on $M$.
\end{theorem}
This statement is proved in
Ref.~\refcite{Ig2} by the fact that
$\mathrm{G}^{1}_{n}=
\mathrm{SL}^{\pm}_{n}\cdot
\mathrm{CO}_{q,n\textrm{-}q}$
and $\mathrm{O}_{q,n\textrm{-}q}=
\mathrm{SL}^{\pm}_{n}\cap
\mathrm{CO}_{q,n\textrm{-}q}$
imply that volume
and conformal $G$-structures
intersect in
$\mathrm{O}_{q,n\textrm{-}q}$-structures.

\section{Second order $G$-structures}\label{SeProl}
A \emph{symmetric linear connection} (SLC)
on $M$ is a distribution on $LM$
of torsion-free transversal $n$-subspaces,
which is invariant by the action of
$\mathrm{G}^{1}_{n}$. Identifying
$\mathrm{G}^{1}_{n}
\simeq \mathrm{G}^{1}_{n}
\rtimes 0 \subset
\mathrm{G}^{2}_{n}$,
we can define an SLC on $M$ as a
second order
$\mathrm{G}^{1}_{n}$-structure.
From bundle theory,
as in Th.\ref{ThVol},
every SLC $\nabla $
is a section of the
associated bundle
to $\mathcal{F}^{2}M$
and the action of
$\mathrm{G}^{ 2}_{n}$ on
$\mathrm{G}^{ 2}_{n}/
\mathrm{G}^{ 1}_{n}
\simeq \mathrm{S}^{2}_{n}$,
and corresponds to an
equivariant function
$f^{\nabla }\colon
\mathcal{F}^{2}M\rightarrow
\mathrm{S}^{2}_{n}$, verifying
$f^{\nabla }(z(a,s))=a^{-1}
f^{\nabla }(z)(a,a)+s$,
for $z\in \mathcal{F}^{2}M$, $a\in
\mathrm{G}^{1}_{n}$,
$s\in \mathrm{S}^{2}_{n}$.

Let $P$ be a first order
$G$-structure on $M$; a
\emph{symmetric connection
on $P$} is a distribution on $P$ of
torsion-free transversal
$n$-subspaces, which is invariant
by the action of $G$,
thereby producing a
second order $G$-structure, whose
$G^{1}_{n}$-extension
is an SLC on $M$.
\textbf{On the other hand, a second order $G$-structure projects onto a first order $G$-structure, say $P$, and by extension it determines an SLC, but it is not, in general, a symmetric connection on $P$.}

Noteworthy examples of this are:
\textbf{i)} A pseudo-Riemannian metric
and its \emph{Levi-Civita connection}
\textbf{determine} a second order
$\mathrm{O}_{q,n\textrm{-}q}$-structure.
\textbf{ii)} An \emph{equiaffine
structure} on $M$ is a SLC
with a parallel volume;
it \textbf{determines} a second order
$\mathrm{SL}^{\pm}_{n}$-structure.
\textbf{iii)} A \emph{Weyl
structure} is a conformal structure
with a SLC compatible;
it \textbf{determines} a second order
$\mathrm{CO}_{q,n\textrm{-}q}$-structure.

The following result
is an important theorem,
arisen from the
Weyl's \emph{`Raumproblem'},
studied by Cartan
and others. The theorem is
proved in Ref.~\refcite{KoNa1}, with
a correction revealed in
Ref.~\refcite{Urba1}.

\begin{theorem}
Let $G$ be a subgroup
of $G^{1}_{n}$,
with $n\geq 3$.
Any first order
$G$-structure admits
a symmetric connection
if and only if
$\mathfrak{g}$
is one of these:
$\mathfrak{sl}(n,\mathbb{R})$,
$\mathfrak{o}_{q,n\textrm{-}q}$,
$\mathfrak{co}_{q,n\textrm{-}q}$,
$\mathfrak{gl}_{n,W}$
(algebra of endomorphisms
with an invariant
1-dimensional subspace $W$),
$\mathfrak{gl}_{n,W,c}$
(certain subalgebra
of the last one, for each
$c\in\mathbb{R}$)
or, for $n=4$,
$\mathfrak{csp}(2,\mathbb{R})$.
\end{theorem}

A $G$-structure $P$ can admit many
torsion-free transversal $n$-subspaces
in $T_{l}P$, for every $l\in P$.
We have the following result (see
\cite[p.150-155]{IgTes}):

\begin{theorem}
If a first order $G$-structure
$P$ admits a symmetric
connection, the set
$P^{2}$
of 2-frames corresponding
with torsion-free
transversal $n$-subspaces
included in $TP$
is a second order $G\rtimes
\mathfrak{g}_{1}$-structure.
\end{theorem}

We name $P^{2}$
the \emph{(first holonomic)
prolongation of} $P$.
\textbf{But} a second order $G\rtimes
\mathfrak{g}_{1}$-structure
is \textbf{not necessarily} the prolongation of
a $G$-structure.

\textbf{Following with the above-mentioned examples:
\textbf{i)} The first holonomic prolongation 
of a pseudo-Riemannian metric
is a second order }
$\mathrm{O}_{q,n\textrm{-}q}$\textbf{-structure,
whose $G^{1}_{n}$-extension is 
the \emph{Levi-Civita connection}.
\textbf{ii)} An \emph{equiaffine
structure} on $M$ can be characterized
as a reduction of the first 
holonomic prolongation of a volume on $M$ to }
$\mathrm{SL}^{\pm}_{n}\subset G^{2}_{n}$.
\textbf{iii)} \textbf{A \emph{Weyl
structure}  can be characterized
as a reduction of the first
holonomic prolongation of 
a conformal structure to a second order}
$\mathrm{CO}_{q,n\textrm{-}q}$\textbf{-structure.}

A \emph{(differential)
projective structure}
is a set of SLCs
which have the same geodesics
up to reparametrizations;
we can define it
as a second order
\emph{$\mathrm{G}^{1}_{n}\rtimes
\mathfrak{p}$-structure} with
$\mathfrak{p}:=\{
s\in \mathrm{S}^{2}_{n}\colon
s^{i}_{jk}=\delta^{i}_{j}\mu _{k}+\mu _{j}
\delta^{i}_{k},\ \mu =(\mu _{i})\in
\mathbb{R}^{n\,\ast }\} \simeq
\mathbb{R}^{n\,\ast }$.
Considering geometrical structures
on second order, with
the same techniques than
in the Th.\ref{ThMetric}
for the first order (see
Ref.~\refcite{Ig2}), we obtain:
\begin{theorem}
A projective structure
and a volume on $M$
give an SLC belonging to the former
and making the volume parallel.
\end{theorem}

Hence, a volume select
a class of affine parametrizations
for the paths of a
projective structure.
Contrarily, a projective structure
and the prolongation of a conformal
structure not always intersect;
if they intersect we get a Weyl
structure.

\section{Concluding remarks}
The study of integrability conditions
of higher order, like curvatures,
with respect to the interrelations of
these $G$-structures
is probably the natural next step
following this work.

The geometrical structures
described herein can be considered
components of the
space-time geometry \cite{Lisboa}.
In this context, the named
\emph{causal set theory}
make a conceptual separation
between volume and conformal
structures, and Stachel \cite{Stac1}
proposes an approach,
similar to the metric-affine variational
principle, using conformal and
projective structures
as independent variables.
In this line of thought,
and from the above results,
I suggest considering the
volume on space-time
as a set of independent
dynamical variables
to make a variational analysis.

\section*{Acknowledgments}
This work has been partially
supported by the \emph{Andalusian
Government} P.A.I.: FQM-324.

\end{document}